\documentclass[11pt]{amsart}
\usepackage{amssymb,amsmath,amsgen,amsxtra,amsfonts,a4wide} 
\usepackage{dsfont,times}

\begin{document}
%
%
%
\theoremstyle{definition}
\newtheorem{Definition}{Definition}[section]
\newtheorem*{Definitionx}{Definition}
\newtheorem{Convention}{Definition}[section]
\newtheorem{Construction}{Construction}[section]
\newtheorem{Example}[Definition]{Example}
\newtheorem{Examples}[Definition]{Examples}
\newtheorem{Remark}[Definition]{Remark}
\newtheorem{Remarks}[Definition]{Remarks}
\newtheorem{Caution}[Definition]{Caution}
\newtheorem{Conjecture}[Definition]{Conjecture}
\newtheorem{Question}[Definition]{Question}
\newtheorem{Questions}[Definition]{Questions}
\theoremstyle{plain}
\newtheorem{Theorem}[Definition]{Theorem}
\newtheorem*{Theoremx}{Theorem}
\newtheorem{Proposition}[Definition]{Proposition}
\newtheorem*{Propositionx}{Proposition}
\newtheorem{Lemma}[Definition]{Lemma}
\newtheorem{Corollary}[Definition]{Corollary}
\newtheorem{Fact}[Definition]{Fact}
\newtheorem{Facts}[Definition]{Facts}
\newtheoremstyle{voiditstyle}{3pt}{3pt}{\itshape}{\parindent}%
{\bfseries}{.}{ }{\thmnote{#3}}%
\theoremstyle{voiditstyle}
\newtheorem*{VoidItalic}{}
\newtheoremstyle{voidromstyle}{3pt}{3pt}{\rm}{\parindent}%
{\bfseries}{.}{ }{\thmnote{#3}}%
\theoremstyle{voidromstyle}
\newtheorem*{VoidRoman}{}

%
\newcommand{\prf}{\par\noindent{\sc Proof.}\quad}
\newcommand{\blowup}{\rule[-3mm]{0mm}{0mm}}
\newcommand{\cal}{\mathcal}
\newcommand{\Aff}{{\mathds{A}}}
\newcommand{\BB}{{\mathds{B}}}
\newcommand{\CC}{{\mathds{C}}}
\newcommand{\FF}{{\mathds{F}}}
\newcommand{\GG}{{\mathds{G}}}
\newcommand{\HH}{{\mathds{H}}}
\newcommand{\NN}{{\mathds{N}}}
\newcommand{\ZZ}{{\mathds{Z}}}
\newcommand{\PP}{{\mathds{P}}}
\newcommand{\QQ}{{\mathds{Q}}}
\newcommand{\RR}{{\mathds{R}}}
\newcommand{\Liea}{{\mathfrak a}}
\newcommand{\Lieb}{{\mathfrak b}}
\newcommand{\Lieg}{{\mathfrak g}}
\newcommand{\Liem}{{\mathfrak m}}
\newcommand{\ideala}{{\mathfrak a}}
\newcommand{\idealb}{{\mathfrak b}}
\newcommand{\idealg}{{\mathfrak g}}
\newcommand{\idealm}{{\mathfrak m}}
\newcommand{\idealp}{{\mathfrak p}}
\newcommand{\idealq}{{\mathfrak q}}
\newcommand{\idealI}{{\cal I}}
\newcommand{\lin}{\sim}
\newcommand{\num}{\equiv}
\newcommand{\dual}{\ast}
\newcommand{\iso}{\cong}
\newcommand{\homeo}{\approx}
\newcommand{\mm}{{\mathfrak m}}
\newcommand{\pp}{{\mathfrak p}}
\newcommand{\qq}{{\mathfrak q}}
\newcommand{\rr}{{\mathfrak r}}
\newcommand{\pP}{{\mathfrak P}}
\newcommand{\qQ}{{\mathfrak Q}}
\newcommand{\rR}{{\mathfrak R}}
%
%
\newcommand{\dq}{{``}}
\newcommand{\OO}{{\cal O}}
\newcommand{\into}{{\hookrightarrow}}
\newcommand{\onto}{{\twoheadrightarrow}}
\newcommand{\Spec}{{\rm Spec}\:}
\newcommand{\BigSpec}{{\rm\bf Spec}\:}
\newcommand{\Proj}{{\rm Proj}\:}
\newcommand{\Pic}{{\rm Pic }}
\newcommand{\Br}{{\rm Br}}
\newcommand{\NS}{{\rm NS}}
\newcommand{\chit}{\chi_{\rm top}}
\newcommand{\KanDiv}{{\cal K}}
\newcommand{\perdef}{{\stackrel{{\rm def}}{=}}}
\newcommand{\Cycl}[1]{{\ZZ/{#1}\ZZ}}
\newcommand{\Sym}{{\mathfrak S}}
\newcommand{\Xcan}{X_{{\rm can}}}
\newcommand{\Ycan}{Y_{{\rm can}}}
\newcommand{\ab}{{\rm ab}}
\newcommand{\Aut}{{\rm Aut}}
\newcommand{\Hom}{{\rm Hom}}
\newcommand{\Supp}{{\rm Supp}}
\newcommand{\ord}{{\rm ord}}
\newcommand{\divisor}{{\rm div}}
\newcommand{\Alb}{{\rm Alb}}
\newcommand{\Jac}{{\rm Jac}}
\newcommand{\piet}{{\pi_1^{\rm \acute{e}t}}}
\newcommand{\Het}[1]{{H_{\rm \acute{e}t}^{{#1}}}}
\newcommand{\Hcris}[1]{{H_{\rm cris}^{{#1}}}}
\newcommand{\HdR}[1]{{H_{\rm dR}^{{#1}}}}
\newcommand{\hdR}[1]{{h_{\rm dR}^{{#1}}}}
\newcommand{\defin}[1]{{\bf #1}}

\title[Non-Classical Godeaux Surfaces]{Non-Classical Godeaux Surfaces}
\author{Christian Liedtke}
\address{Mathematisches Institut, Heinrich-Heine-Universit\"at, 40225
  D\"usseldorf, Germany}
\email{liedtke@math.uni-duesseldorf.de}
\date{April 21, 2008. {\it revised:} August 4, 2008}

\begin{abstract}
 A non-classical Godeaux surface is a minimal surface of general type 
 with $\chi=K^2=1$ but with $h^{01}\neq0$.
 We prove that such surfaces fulfill $h^{01}=1$
 and they can exist only over fields of positive characteristic 
 at most $5$.
 Like non-classical Enriques surfaces they fall into two classes:
 the singular and the supersingular ones.
 We give a complete classification in characteristic $5$ and 
 compute their Hodge-, Hodge--Witt- and crystalline cohomology
 (including torsion).
 Finally, we give an example of a supersingular Godeaux surface
 in characteristic $5$.
\end{abstract}
\setcounter{tocdepth}{1}
\maketitle
\tableofcontents
\section*{Introduction}

Among minimal surfaces of general type, the ones with the smallest invariants
possible are those with $K^2=1$ and $h^{01}=p_g=0$.
The first example of such a surface has been constructed 
by Godeaux \cite{godeaux} in 1931.
Since then, there is an extensive search for such surfaces and up to 
today there is still no complete classification of them.
We refer to \cite[Chapter VII.10]{bhpv} for details and further references.

\begin{Definitionx}
  A minimal surface of general type with $\chi(\OO_X)=K_X^2=1$
  is called a {\em numerical Godeaux surface}.
  A numerical Godeaux surface with $h^{01}(X)\neq0$
  is called {\em non-classical}.
\end{Definitionx}

\begin{Propositionx}
 The first Betti number of a numerical Godeaux surface is zero.
 Non-classical Godeaux surfaces fulfill $h^{01}=p_g=1$.
\end{Propositionx}

This means that the Picard scheme of a non-classical Godeaux surface is not
reduced, which can happen over fields of positive characteristic
only.
In particular, a numerical Godeaux surface 
in characteristic zero fulfills 
$h^{01}=p_g=0$.
On the other hand, Miranda \cite{mir} has given an
example of a non-classical Godeaux surface in 
characteristic $5$.
Our first main result is the following

\begin{Theoremx}
 Non-classical Godeaux surfaces can exist in characteristic
 $2\leq p\leq 5$ only. 
\end{Theoremx}

This should be compared with Enriques surfaces, i.e.,
minimal surfaces of Kodaira dimension zero with $\chi(\OO_X)=1$.
The classical ones have $h^{01}=p_g=0$, whereas the non-classical
ones fulfill $h^{01}=p_g=1$.
Non-classical Enriques surfaces exist in characteristic $2$ only
\cite[Theorem 5]{bm2}.
\medskip

Like Enriques surfaces, all global $1$-forms on a numerical 
Godeaux surface are $d$-closed.
Another feature of Enriques surfaces is that their second
crystalline or \'etale cohomology is spanned by algebraic cycles, i.e.,
that they are supersingular in the sense of Shioda.
This is also true for numerical Godeaux surfaces
if the surface lifts to characteristic zero or if a conjecture
of Artin and Mazur holds true, cf. Section \ref{supersingularsection}.
\medskip

Since $H^1(\OO_X)$ of a non-classical Godeaux surface is one-dimensional,
the action of Frobenius on it is either zero or bijective.
In analogy with Enriques surfaces \cite[Section 3]{bm3} we define

\begin{Definitionx}
  A numerical Godeaux surface is called
  \begin{center}
  \begin{tabular}{ll}
 ×  {\em classical} & if $p_g=h^{01}=0$,\\
 ×  {\em singular} & if $p_g=h^{01}=1$ and $F$ is bijective on $H^1(\OO_X)$, and\\
 ×  {\em supersingular} & if $p_g=h^{01}=1$ and $F$ is zero on $H^1(\OO_X)$.
  \end{tabular}
  \end{center}
\end{Definitionx}

Now we specialise to characteristic $5$, the largest characteristic possible for
non-classical Godeaux surfaces.
As already mentioned, the main feature is the non-reducedness of 
their Picard schemes.
More precisely, we show

\begin{Theoremx}
   Let $X$ be a non-classical Godeaux surface in characteristic $5$.
   Then
  \begin{center}
    \begin{tabular}{lcl}
     $X$ singular &implies& $\Pic^\tau(X)\,=\,\Pic^0(X)\,\iso\,\mu_5$\,, and\\
     $X$ supersingular &implies& $\Pic^\tau(X)\,=\,\Pic^0(X)\,\iso\,\alpha_5\,.$
    \end{tabular}
  \end{center}
\end{Theoremx}

This result is the key to determining the Hodge-, the Hodge--Witt- and the 
crystalline cohomology groups (including torsion) of non-classical Godeaux surfaces 
in characteristic $5$.
Also the degeneration behaviour of the Fr\"olicher- and the slope spectral sequence are
answered.
A picture emerges that is very similar to non-classical Enriques surfaces.
We refer to Section \ref{characteristic5} for precise statements.

Next, we prove that every non-classical Godeaux surface
in characteristic $5$ arises as the
quotient of a quintic surface in $\PP^3$ by a linear and fixed-point free action of
$\ZZ/5\ZZ$, resp. $\alpha_5$, i.e., every such surface is obtained
by a characteristic $5$-version of Godeaux's original construction.
This result is analogous to Reid's classification 
\cite[Section 1]{reid} of numerical Godeaux 
surfaces with $5$-torsion over the complex numbers.

\begin{Theoremx}
 Let $\Xcan$ be the canonical model of a non-classical Godeaux surface in characteristic $5$.
 Then its associated $\ZZ/5\ZZ$-torsor (in the singular case), resp.
 $\alpha_5$-torsor (in the supersingular case)
 is isomorphic to a quintic surface in $\PP^3$.
 Moreover, the $\ZZ/5\ZZ$-, resp. $\alpha_5$-action on this quintic extends
 to a linear action on the ambient $\PP^3$.
\end{Theoremx}

Finally, it remains to establish existence of these surfaces.
Miranda's surface \cite{mir} mentioned above 
is a singular Godeaux surface in characteristic $5$
and we end this article by proving

\begin{Theoremx}
 There do exist supersingular Godeaux surfaces in characteristic $5$.
\end{Theoremx}

\begin{VoidRoman}[Acknowledgements]
 I thank Stefan~Schr\"oer for interesting discussions and help.
 Also, I thank the referee for pointing out a couple of inaccuracies
 and a mistake in the first version of this article.
\end{VoidRoman}

\section{Numerical Godeaux surfaces}

A minimal surface $X$ of general type 
defined over an arbitrary algebraically closed field
fulfills $K_X^2\geq1$.
We first recall a couple of general facts about
surfaces for which the equality $K_X^2=1$ holds.
In particular, we will see that $\chi(\OO_X)\geq1$ holds true.
Hence surfaces with $K_X^2=\chi(\OO_X)=1$,
the so-called {\em numerical Godeaux surfaces},
have the lowest invariants possible among all
minimal surfaces of general type.
Although the first Betti number of such a surface is always
zero, over fields of positive characteristic we only
have $h^{01}\leq1$, which leads to the notion of
{\em non-classical Godeaux surfaces}.
\medskip

We start with some general facts that are well known over the
complex numbers, cf. Sections 10 and 11 of \cite{bom}.
However, a little care is needed in positive characteristic.

\begin{Proposition}
 \label{general}
 Let $X$ be a minimal surface of general type with $K_X^2=1$.
 Then the following equalities and inequalities hold true:
 $$
  b_1(X)=0,\,\mbox{ \quad }\,
  |\piet(X)|\leq6,\,\mbox{ \quad }\,
  1\leq\chi(\OO_X)\leq3,\,\mbox{ \quad }\,
  p_g(X)\leq2,\,\mbox{ \quad and, \quad }\,
  h^{01}(X)\leq1\,.
 $$
 In particular, if $h^{01}(X)\neq0$ then $X$ has a non-reduced
 Picard scheme, which can happen in positive characteristic only.
\end{Proposition}

\prf
From Noether's inequality $K^2\geq2p_g-4$,
which also holds in positive characteristic by
\cite[Theorem 2.1]{lie}, we obtain $p_g\leq2$
and then $\chi\leq3$.
Then, $K^2+\chi\geq2$ yields
$\chi\geq1$, cf. \cite[Corollary II.1.8]{ek}.
Using the inequalities of \cite[Corollary II.1.8]{ek} again, 
we see that $h^{01}\leq1$ if $p_g=0$ or $p_g=2$.
In case $p_g=1$, we also must have $h^{01}\leq1$ since
$\chi\geq1$.

Let $\tilde{X}$ be an \'etale cover of $X$ of degree $n$.
Then $\chi(\OO_{\tilde{X}})=n\chi(\OO_X)\geq n$,
which implies $p_g(\tilde{X})\geq n-1$.
We also have $K_{\tilde{X}}^2=n$ and then Noether's inequality
$K^2\geq2p_g-4$ implies $n\leq6$.
Hence the \'etale fundamental group of $X$ has order at most $6$,
which implies that the first Betti number of $X$ is zero.
\qed\medskip

Since we are interested in surfaces of general type with the lowest
invariants possible, the previous result tells us to look at the 
following class of surfaces:

\begin{Definition}
  A {\em numerical Godeaux surface} is a minimal surface of general
  type with $K_X^2=1$ and $\chi(\OO_X)=1$.
\end{Definition}

If a numerical Godeaux surface fulfills $p_g=h^{01}\neq0$ we are in
positive characteristic and it makes sense to look at the action
$F$ of Frobenius on $H^1(\OO_X)$, which is either zero or bijective
since this space is one-dimensional.
In analogy with Enriques surfaces \cite[Section 3]{bm3} we define:

\begin{Definition}
  A numerical Godeaux surface is called
  \begin{center}
  \begin{tabular}{ll}
 ×  {\em classical} & if $p_g=h^{01}=0$,\\
 ×  {\em singular} & if $p_g=h^{01}=1$ and $F$ is bijective on $H^1(\OO_X)$, and\\
 ×  {\em supersingular} & if $p_g=h^{01}=1$ and $F$ is zero on $H^1(\OO_X)$.
  \end{tabular}
  \end{center}
  We will refer to the last two classes as
  {\em non-classical Godeaux surfaces}.
\end{Definition}

By Proposition \ref{general}, every numerical Godeaux surface in characteristic 
zero is classical.
\medskip

Godeaux's original construction \cite{godeaux} 
of a classical Godeaux surface is 
the quotient of a smooth quintic in $\PP^3$ by a fixed point 
free action of $\mu_5$, which works over
arbitrary fields of characteristic $p\neq5$.
Lang \cite{la} also
constructed classical Godeaux surfaces in characteristic $5$.

\begin{Theorem}[Godeaux, Lang]
  There do exist classical Godeaux surfaces in every characteristic.
\qed
\end{Theorem}

We remark that
the classification of numerical Godeaux surfaces even over
the complex numbers is still incomplete,
and refer to \cite[Chapter VII.10]{bhpv} for details and further
references.

\section{Non-classical Godeaux surfaces}

In the previous section we have seen that classical numerical Godeaux surfaces
exist in every characteristic.
In this section we show that non-classical Godeaux surfaces
can only exist in characteristic at most $5$.
This is rather easy for the singular Godeaux surfaces but much more complicated
for the supersingular Godeaux surfaces.

\begin{Theorem}
  \label{ordinary}
   Singular Godeaux surfaces can exist in characteristic
   $2\leq p\leq5$ only.
   These surfaces are not weakly unirational.
\end{Theorem}

\prf
Let $X$ be a singular Godeaux surface.
Since $H^1(\OO_X)$ is one-dimensional and Frobenius acts injectively, we
obtain an inclusion of the group scheme $\mu_p$ into $\Pic(X)$.
This defines a $\mu_p^D$-torsor over $X$, where $-^D:=\Hom(-,\GG_m)$
denotes Cartier duality.
Since $\mu_p^D\iso\ZZ/p\ZZ$ this torsor is just
an \'etale cover of degree $p$ of $X$.
Hence $p\leq6$ by Proposition \ref{general}.

In particular, the algebraic fundamental group has non-trivial $p$-torsion.
By \cite[Theorem 2.5]{cr}, such a surface cannot be weakly unirational.
\qed\medskip

\begin{Remark}
 Shioda has given examples of classical Godeaux surfaces
 that are weakly unirational in every positive characteristic $p\neq5$
 with $p\not\equiv1\mod5$, cf. \cite[Proposition 5]{sh2}.
\end{Remark}
\medskip

What makes the analysis of supersingular Godeaux surfaces more difficult is
the fact that the trivial Frobenius action on $H^1(\OO_X)$ gives rise to an
$\alpha_p$-torsor over $X$, whose total space need not be a normal surface.
In particular, we cannot apply geometric arguments to this torsor
so easily. 
The main ingredient in the proof of Theorem \ref{ordinary} is Noether's
inequality which we now establish in a singular context:

Let $X$ be a smooth surface
over an algebraically closed field of characteristic $p>0$
and $\pi:Y\to X$ be a non-trivial $\alpha_p$- or $\mu_p$-torsor.
From the discussion \cite[Section I]{ek} it follows that $Y$ is a possibly
non-normal but integral Gorenstein surface.
Since the dualising sheaf $\omega_Y$ of $Y$ is invertible, we define
the self-intersection of the canonical divisor of $Y$ to be 
\begin{equation}
 \label{selfintersection}
  K_Y^2 \,:=\,\chi(\OO_Y)\,-\,2\,\chi(\omega_Y^{-1})\,-\,\chi(\omega_Y^{-2}),
\end{equation}
which coincides with the usual self-intersection in the smooth case.

In \cite[Proposition I.1.7]{ek} the algebra $\pi_\ast\OO_Y$ is shown to be a
successive extension of invertible sheaves.
If $\pi$ is an $\alpha_p$-torsor, then these invertible sheaves are
isomorphic to $\OO_X$.
Also, $\omega_Y\iso\pi^\ast(\omega_X)$ in this case.
If $\pi$ is a $\mu_p$-torsor, then $\pi_\ast\OO_Y$ is a
successive extension of ${\cal L}^{\otimes i}$, $i=0,...,p-1$
for a suitable invertible sheaf $\cal L$ with 
${\cal L}^{\otimes p}\iso\OO_X$.
In this case we have 
$\omega_Y\iso\pi^\ast(\omega_X\otimes{\cal L}^{\otimes(p-1)})$.
It follows that the equalities
\begin{equation}
\label{standard}
\chi(\OO_Y) \,=\,p\,\chi(\OO_X)\,\mbox{ \quad and\quad }
K_Y^2\,=\,p\,K_X^2
\end{equation}
hold true.

What makes Noether's inequality for $Y$ a little bit tricky is
the fact that $Y$ is only an integral Gorenstein scheme, 
so that linear systems and intersection numbers
have to be treated with care.

\begin{Proposition}
 \label{noether}
 Let $X$ be a minimal surface of general type and 
 $\pi:Y\to X$ be a non-trivial $\alpha_p$- or $\mu_p$-torsor.
 Then Noether's inequality
 $$
   K_Y^2\,\geq\,2\,h^0(\omega_Y)\,-\,4
 $$
 holds true.
 If the image of the canonical map is a curve, even 
$K_Y^2\geq2h^0(\omega_Y)-2$ holds true.
\end{Proposition}

\prf
Since $\pi:Y\to X$ is purely inseparable of degree $p$, the 
Frobenius morphism $F$ factors as $\rho: X^{(p)}\to Y$ 
followed by $\pi$.
We already noted that 
$\omega_Y\iso\pi^\ast(\omega_X\otimes{\cal L}^{\otimes(p-1)})$
for a suitable invertible sheaf $\cal L$.
Using the fact that 
${\cal L}^{\otimes p}\iso\OO_X$, we get
$$
\rho^\ast\omega_Y\,\iso\, 
F^\ast(\omega_X\otimes{\cal L}^{\otimes (p-1)})
\,\iso\, 
\omega_{X^{(p)}}^{\otimes p}\,.
$$
We denote by $L$ the linear system $\rho^\ast H^0(\omega_Y)$
of $H^0(\omega_{X^{(p)}}^{\otimes p})$
and consider the composed (rational) morphism
$$
\begin{array}{ccccccccc}
\varphi\,:=\,\varphi_{|L|} &:& X^{(p)} &\stackrel{\rho}{\to}&
Y&\stackrel{\phi_1}{\dashrightarrow}&Z&\subseteq&\PP^{h^0(\omega_Y)-1},
\end{array}
$$
where $Z$ denotes the closure of the image of the (rational) canonical map 
$\phi_1$ of $Y$.
Being a minimal surface of general type we have $K_X^2\geq1$, hence
$K_Y^2\geq p$ and thus we may assume  $h^0(\omega_Y)\geq3$, i.e.,
that $Z$ is a curve or a surface.

Let us first assume that $Z$ is a surface, in which case
$\deg\varphi\geq p$. 
If $\phi_1$ is not birational (or $Z$ is not ruled) then
\cite[Proposition 0.1.3]{ek}, which equally well works for linear
subsystems of complete linear systems, yields 
$h^0(\omega_Y)=\dim L\leq (pK_{X^{(p)}})^2/(2p) + 2=K_Y^2/2+2$ and
we are done.
We may now assume that $\phi_1$ is birational 
(and $Z$ is ruled).
For every section of $\omega_Y$ we obtain a short exact sequence
$$
0\,\to\,\OO_Y\,\to\,\omega_Y\,\to\,\omega_Y|_D\,\to\,0\,,
$$
where $D$ is a Cartier divisor on a Gorenstein scheme.
In particular, $D$ is Gorenstein and the adjunction formula holds.
The divisor $\rho^\ast D$ lies in $\omega_{X^{(p)}}^{\otimes p}$ 
and since 
$h^1(\omega_{X^{(p)}}^{\otimes(-p)})=0$ by \cite[Theorem II.1.7]{ek},
we obtain $h^0(\OO_{\rho^\ast D})=1$.
Thus, $h^0(\OO_D)=1$ and we compute
$h^0(\omega_D)=1-\chi(\OO_D)=1+K_Y^2$.
A generic hyperplane section of $Z$ is an integral curve, and
thus corresponds to a Cartier divisor $D$ on $Y$ which has
a component $D'$, possibly only a Weil divisor, which is an integral
curve (as $\phi_1$ is an isomorphism over an open
and dense subset).
We denote by $V$ the image of
$H^0(\omega_Y)\to H^0(D,\omega_Y|_D)$.
Every non-zero element of $V$ corresponds
to a hyperplane section of $Z$ and is thus non-trivial 
on $D'$.
Since $D'$ is integral, the multiplication
of two non-zero sections of an invertible sheaf is non-zero,
we can apply the Clifford argument and obtain
the estimate
$2h^0(\omega_Y)-4=2\dim V-2\leq h^0(D,\omega_Y|_D^{\otimes2})-1
= h^0(D,\omega_D) - 1 = K_Y^2$.

We may now assume that $Z$ is a curve.
To begin with, we follow the proof of \cite[Theorem VII.3.1]{bhpv}.
We consider the linear subsystem $L:=\rho^\ast H^0(\omega_Y)$
of $H^0(\omega_{X^{(p)}}^{\otimes p})$ and write this pencil
as $|L|=|nF|+V$, where $V$ denotes its fixed part.
Then $h^0(\omega_Y)=\dim L\leq 1+n$.
We compute $K_Y^2=pK_X^2=(nF+V)K_X\geq nFK_X$ and we are done
if $FK_X\geq2$.
In the remaining cases we either have $FK_X=0$ or $FK_X=1$.
If $FK_X=0$ then $F$ is a sum of fundamental cycles and these
cannot form a pencil, i.e., this case does not occur.
It remains to exclude the possibility $FK_X=1$.
Since $F$ moves in a pencil we have $F^2\geq0$ and also $F^2$
must be an odd integer by the adjunction formula, hence
$F^2\geq1$.
The Hodge index theorem yields $F^2K_X^2\leq(FK_X)^2=1$, which
leads to $K_X^2=FK_X=F^2=1$.

It remains to deal with this case:
since $F^2=1$, we have to blow up once $\tilde{X}\to X$ to obtain
a morphism $\tilde{X}\to\PP^1$, which is in fact a fibration.
The $(-1)$-curve $E\subset \tilde{X}$ dominates $\PP^1$, i.e.,
defines a section.
The induced $\alpha_p$-torsor (resp. $\mu_p$-torsor) 
$\tilde{Y}\to\tilde{X}$ restricts to a trivial torsor over $E$
because there are no non-trivial such torsors over a 
rational curve.
In particular, we can extend the section of $\tilde{X}\to\PP^1$
to $\tilde{Y}$.
This means that $\tilde{Y}\to\PP^1$ is equal to its Stein
factorisation, hence that its fibres are geometrically integral.
Hence for every irreducible fibre $F$ of $\tilde{X}\to\PP^1$
its inverse image $\pi^\ast F$ is an integral curve.

We now compute on $X$:
let $F'$ be a general member of $|F|$ and set
$\tilde{F}':=\pi^\ast F'$, which is an integral curve
by what we have just proved.
We consider the short exact sequence
$$
0\,\to\,\omega_Y(-\tilde{F}')\,\to\,\omega_Y\,\to\,
\omega_Y|_{\tilde{F}'}\,\to\,0
$$
Taking cohomology and using Serre duality, we get
$h^1(\omega_Y|_{\tilde{F}'})\neq0$.
Since $\tilde{F}'$ is an integral curve, Clifford's
inequality gives 
$h^0(\omega_Y|_{\tilde{F}'})\leq 1+(\deg\omega_Y|_{\tilde{F}'})/2=1+K_Y^2/2$.
Now, ${\cal B}:=\omega_Y(-\tilde{F}')$ is the pull-back of a numerically trivial
invertible sheaf on $X$, and hence $h^0({\cal B})\leq1$.
From this we get Noether's inequality.

In case $h^0({\cal B})=1$, the curve $\tilde{F}'$ is a reduced section of 
$\omega_Y$ and so the map associated to $\omega_Y$ is not composed with
a pencil.
Hence if the image of the canonical map is a curve, we have
$h^0({\cal B)}=0$ and the stronger inequality $h^0(\omega_Y)\leq 1+K_Y^2/2$
holds true.
\qed\medskip

\begin{Theorem}
   \label{supersingular}
   Supersingular Godeaux surfaces can exist in characteristic
   $2\leq p\leq5$ only.
\end{Theorem}

\prf
Let $X$ be a supersingular Godeaux surface.
By definition, the Frobenius-action on the one-dimensional
vector space $H^1(\OO_X)$ is trivial.
This action gives rise to a non-trivial
$\alpha_p$-torsor $\pi:Y\to X$.
By (\ref{standard}) we have $K_Y^2=pK_X^2=p$ and 
$\chi(\OO_Y)=p\chi(\OO_X)=p$.
The latter equality gives $h^0(\omega_Y)\geq p-1$.
By Proposition \ref{noether},
the inequality $p\geq2(p-1)-4$ has to be fulfilled
and we obtain $p\leq5$.
\qed\medskip

\begin{Remark}
 \label{singularremark}
 The same argument as in \cite[Corollary 1 of Theorem 1]{la}
 shows that the total space of the $\alpha_p$-torsor 
 associated to a supersingular Godeaux surface is never
 a smooth surface.
\end{Remark}

\begin{Corollary}
   Non-classical Godeaux surfaces can exist in characteristic
   $2\leq p\leq5$ only.\qed
\end{Corollary}

\section{Supersingularity and global $1$-forms}
\label{supersingularsection}

In this section we show that numerical Godeaux surfaces that
are liftable to characteristic zero are supersingular in the
sense of Shioda.
Without the lifting assumption
we can show this result only modulo a conjecture
of Artin and Mazur.
We finish by showing that global $1$-forms on a numerical
Godeaux surface are always $d$-closed.
\medskip

By a theorem of Igusa, the rank of the N\'eron--Severi group of a smooth 
projective surface is less or equal to its second Betti number.
This led  Shioda \cite{sh} to introduce 
the following notion:

\begin{Definition}
  A smooth and projective surface is called 
  {\em supersingular in the sense of Shioda} if
  the rank of its N\'eron--Severi group is equal to its second
  Betti number.
\end{Definition}

Over fields of characteristic zero, a smooth projective
surface is supersingular in the 
sense of Shioda if and only if $p_g=0$ by \cite[Section 2]{sh2}.

\begin{Proposition}
 If a numerical Godeaux surfaces
 is defined over a field of characteristic zero or if it
 lifts to characteristic zero 
 then it is supersingular in the sense of Shioda.
\end{Proposition}

\prf
In characteristic zero we have $p_g=0$ and the 
result follows from \cite[Section 2]{sh2}.

If a lifting of $X$ 
to characteristic zero exists, then the lifted
surface $X^\ast$ also has $K^2=\chi=1$.
Hence $p_g(X^\ast)=0$ by Proposition \ref{general} 
and thus $X^\ast$ is supersingular in the sense of Shioda.
We denote by $\rho$ the rank of the N\'eron--Severi group
and conclude as in
\cite[Lemma 1]{sh2}
$$
b_2(X^\ast)\,=\,b_2(X)\,\geq\,\rho(X)\,\geq\,\rho(X^\ast)\,=\,b_2(X^\ast)\,.
$$
Hence $\rho(X)=b_2(X)$, i.e., $X$ is supersingular in the sense
of Shioda.
\qed\medskip

If a surface over an algebraically closed field of positive
characteristic is supersingular in the sense of Shioda,
then the Frobenius-action on its second crystalline cohomology
modulo torsion is of slope one, 
cf. \cite[Proposition II.5.12]{ill}. 

Artin and Mazur conjectured that the converse is also true,
cf. \cite[Remarque II.5.13]{ill}.
For example, this is true for Enriques surfaces \cite[Theorem 4]{bm3}.
Modulo this conjecture, numerical Godeaux surfaces
are supersingular in the sense of Shioda:

\begin{Proposition}
 \label{slopeone}
  In positive characteristic, the second 
  crystalline cohomology modulo torsion 
  of a numerical Godeaux surface is of slope one.
\end{Proposition}

\prf
As in the proof of \cite[Proposition II.7.3.2]{ill} one
first shows that the Witt vector cohomology $H^2(W\OO_X)$
is equal to its $V$-torsion.
Now, this $V$-torsion is a $W$-module of finite length since
it can be related to the Dieudonn\'e-module
of some finite flat group scheme, cf.
\cite[Remarque II.6.4]{ill}.
Hence $H^2(W\OO_X)$ is a $W$-module of finite length, which
implies that $\Hcris{2}(X/W)$ modulo torsion is of slope one,
cf. \cite[Formula II.(3.5.4)]{ill}.
\qed\medskip

As in \cite[Corollaire II.7.3.3]{ill} we obtain

\begin{Proposition}
  \label{dclosed}
  The slope spectral sequence of a numerical Godeaux
  surface degenerates at $E_1$-level.
  In particular,
  every global $1$-form on a numerical Godeaux surface is $d$-closed.
  \qed
\end{Proposition}

Clearly, there are no non-trivial global $1$-forms on a numerical 
Godeaux surface over the complex numbers since $b_1=0$.
On the other hand, supersingular Godeaux surfaces and
Lang's classical Godeaux surfaces in characteristic $5$ 
possess $\alpha_p$-, resp. $\mu_p$-torsors above them, which 
implies that they have non-trivial global $1$-forms,
cf. \cite[Proposition I.0.1.11]{cd} and Proposition \ref{e1terms}
below.

\section{Picard scheme and Hodge invariants}
\label{characteristic5}

We have seen that non-classical Godeaux surfaces fall into two
classes: the singular and the supersingular ones.
Also we have seen that these surfaces can exist in
characteristic $2\leq p\leq5$ only.
In this section we specialise to characteristic $5$.
We will determine their $\Pic^\tau$'s and their Hodge invariants 
as well as their crystalline cohomology.
As we will see, these surfaces behave in many respects like
non-classical Enriques surfaces.
\medskip

\begin{Theorem}
 \label{picard}
  Let $X$ be a non-classical Godeaux surface in characteristic $5$.
  Then
  \begin{center}
    \begin{tabular}{lcl}
     $X$ singular &implies& $\Pic^\tau(X)\,=\,\Pic^0(X)\,\iso\,\mu_5$\,, and\\
     $X$ supersingular &implies& $\Pic^\tau(X)\,=\,\Pic^0(X)\,\iso\,\alpha_5\,.$
    \end{tabular}
  \end{center}
 Singular Godeaux surfaces fulfill $\piet(X)\iso(\ZZ/5\ZZ)$ and
 supersingular Godeaux surfaces are algebraically simply connected.
\end{Theorem}

\prf
First, we determine $\Pic^0$.
Assume that $X$ is supersingular.
Then we have an inclusion of $\alpha_5$ into $G:=\Pic^0(X)$
giving rise to an $\alpha_5^D$-torsor $\pi:Y\to X$.
If $G':=G/\alpha_5$ were non-trivial then we would find 
an embedding of $\alpha_5$ or $\mu_5$ into $G'$, i.e., there would be
$\alpha_5^D$- or $\mu_5^D$-torsors above $Y$.
In particular, we would have $h^1(\OO_Y)\neq0$.
We have $K_Y^2=\chi(\OO_Y)=p=5$ by (\ref{standard}),
which yields $h^0(\omega_Y)\geq4$.
Then Proposition \ref{noether} forces $h^0(\omega_Y)\leq4$,
whence $h^0(\omega_Y)=4$ and hence $h^1(\OO_Y)=0$.
This contradiction shows that $G'$ is trivial, i.e.,
$\Pic^0\iso\alpha_5$.
If $X$ is singular, we only have to replace $\alpha_5$ by $\mu_5$
in the previous discussion.

Next, we show that there is no \'etale $5$-torsion in $\Pic^\tau$.
Otherwise, a non-trivial invertible sheaf $\cal L$ with 
$\cal L^{\otimes 5}\iso\OO_X$ would give rise to a
$\mu_5$-torsor $\pi:Z\to X$.
As before we conclude $h^1(\OO_Z)=0$.
If $X$ is supersingular, we have an inclusion 
of $\alpha_5\times(\ZZ/5\ZZ)$ into
$\Pic^\tau$, which yields a non-trivial
$\alpha_5^D$-torsor over $Z$.
But then $h^1(\OO_Z)\neq0$, which we have just excluded.
In case $X$ is singular, we argue with $\mu_5\times(\ZZ/5\ZZ)$
and conclude as before.

As there is no \'etale $5$-torsion,
$Q:=\Pic^\tau(X)/\Pic^0(X)$ is an \'etale
group scheme of order $\ell(Q)$ prime to $5$.
It remains to show that $Q$ is trivial.

If $X$ is singular, then $\Pic^\tau$ gives rise to an 
\'etale cover of degree $5\cdot\ell(Q)$.
Proposition \ref{general} gives $5\cdot\ell(G)\leq6$, 
and hence $Q$ is trivial.
Moreover, since we already have an inclusion of $\ZZ/5\ZZ$
into $\piet(X)$, Proposition \ref{general} 
tells us that $\piet(X)$ is in fact isomorphic
to $\ZZ/5\ZZ$.

We will assume from now on that $X$ is supersingular.
We know that $\ell(Q)\leq6$ and $\ell(Q)\neq5$.
If $\ell(Q)=6$ then $Q\iso\ZZ/6\ZZ$.
Hence if we show that there are no elements of order
$2$ and $3$ in $\Pic^\tau$, the case $\ell(Q)=6$
is excluded at the same time.

A hypothetical element of order $2$ in $\Pic^\tau(X)$ would give rise
to an \'etale cover $\varpi:Z\to X$ of degree $2$
and there would still exist a non-trivial 
$\alpha_5$-torsor $Y\to Z$.
We have $\chi(\OO_Y)=K_Y^2=2\cdot 5=10$ and in
particular $h^0(\omega_Y)\geq 9$.
Applying the inequality of Proposition \ref{noether} to the
$\alpha_5$-torsor $Y\to Z$ we obtain a contradiction.
Hence there is no non-trivial $2$-torsion in $\Pic^\tau(X)$.
The case of $3$-torsion is excluded similarly.

From Proposition \ref{general} we know that $\piet(X)$
is of order at most $6$.
Also this group is of order prime to $5$ as $X$ is
supersingular.
If it were non-trivial then it would have to have
cyclic quotients of order $2$ or $3$.
These quotients would give us elements of order
$2$ or $3$ in $\Pic^\tau$ which we have
just excluded.
Hence $\piet(X)$ is trivial.
\qed\medskip

Theorem \ref{picard} together with Proposition \ref{dclosed} allows us
to compute the Hodge-, the deRham- and the crystalline cohomology groups.
The computations for Enriques surfaces presented in
\cite[Section II.7.3]{ill} carry over literally, which is
why we omit the proof and only state the result:

\begin{Proposition}
  Let $X$ be a non-classical Godeaux surface in characteristic $5$.
  Then
  $$
   \hdR{1}(X/k)\,=\,1,\mbox{ \quad } 
   h^1(\OO_X)\,=\,1\,\mbox{ \quad and \quad } 
   h^0(\Omega_X^1)\,=\,\left\{\begin{tabular}{ll}
                               0 & \mbox{ if $X$ is singular,}\\
                               1 & \mbox{ if $X$ is supersingular.}
                              \end{tabular}
 \right.
  $$  
  The Fr\"olicher spectral sequence 
  degenerates at $E_1$-level if and only if the surface is singular.
  The crystalline cohomology groups are given by
  $$
   \begin{array}{ccccccccc}
     \blowup
     \Hcris{0}(X/W) &\mbox{ }& \Hcris{1}(X/W) &\mbox{ }& \Hcris{2}(X/W) 
     &\mbox{ }& \Hcris{3}(X/W) &\mbox{ }& \Hcris{4}(X/W)\\
     W && 0 && W^9\oplus k && k && W
   \end{array}
  $$
  In any case,
  the slope spectral sequence degenerates at $E_1$-level.\qed
\end{Proposition}

Also, \cite[Proposition II.7.3.6]{ill} and \cite[Proposition II.7.3.8]{ill} carry 
over literally.
We only remark that we do not know whether $\rho=b_2$ holds true for non-classical
Godeaux surfaces but when this result is used in the proof of 
\cite[Proposition II.7.3.6]{ill}, we only need that $\Hcris{2}(X/W)$ modulo torsion
is of slope one, which we have established for numerical Godeaux surfaces
in Proposition \ref{slopeone}.

\begin{Proposition}
 \label{e1terms}
 Let $X$ be a non-classical Godeaux surface in characteristic $5$.
 Then the $E_1$-terms of the Fr\"olicher- and the slope spectral
 spectral sequence are given by
 \begin{center}
  \begin{tabular}{lcccccccc}
   & \qquad & Fr\"olicher spectral sequence& \qquad &slope spectral sequence\\
   \\
   singular && 
             $\begin{array}{ccc}
              k & 0 & k \\
              k & k^9 & k \\
              k & 0 & k
              \end{array}$
           && $\begin{array}{ccc}
              k & 0   & W \\ 
              0 & W^9 & k \\ 
              W & 0 & 0 
              \end{array}$ \\ 
   \\
   supersingular && $\begin{array}{ccc}
                      k & k & k\\
                      k & k^{11} & k\\
                      k & k & k
                    \end{array}$
            && $\begin{array}{ccc}
                      k & k & W \\
                      0 & W^9 & 0 \\
                      W & 0 & 0
                 \end{array}$
  \end{tabular}
 \end{center}
  Finally, we have $\HdR{2}(X/k)=k^{11}$ and $\HdR{3}(X/k)=k$.\qed
\end{Proposition}

\section{Classification and examples}

In this section we show that non-classical Godeaux surfaces in 
characteristic $5$ are obtained as quotients of quintic hypersurfaces
in $\PP^3$ by a linear fixed point free $\ZZ/5\ZZ$-action (singular Godeaux
surfaces), resp. $\alpha_5$-action (supersingular Godeaux surfaces).
This is analogous to Reid's classification 
\cite[Section 1]{reid} of numerical Godeaux surfaces
with $5$-torsion over the complex numbers.
As Miranda \cite{mir} has already shown existence of singular
Godeaux surfaces, we finish this article by establishing the existence
of supersingular Godeaux surfaces.

We start with a classification result similar to the one 
of \cite[Section 2]{la} for classical Godeaux surfaces in 
characteristic $5$.
The main difficulty is the analysis of supersingular Godeaux
surfaces, because the total space of the associated
$\alpha_5$-torsor is a possibly non-normal surface.
If this total space were smooth (which is the case for
singular Godeaux surfaces) we could simply argue along
the lines of \cite[Lemma2]{horquint}
and \cite[Theorem 1]{horquint}.

\begin{Theorem}
 \label{classification}
 Let $X$ be a non-classical Godeaux surface in characteristic $5$
 and $\Xcan$ its canonical model.
 We denote by $Y\to X$ the associated $\ZZ/5\ZZ$-torsor (in case $X$
 is singular), resp. $\alpha_5$-torsor (in case $X$ is supersingular).
 Then there exists a Cartesian diagram 
 $$
 \begin{array}{ccccc}
   Y &\to& Y' &\stackrel{\phi_1}{\to}& \PP^3 \\
   \downarrow&&\downarrow\\
   X &\to&\Xcan
 \end{array}
 $$
 where $Y'$ is the canonical model of $Y$.
 The canonical map of $Y$ is a birational morphism onto a quintic surface in
 $\PP^3$.
 Moreover, it factors over $Y'$ and induces an isomorphism of 
 $Y'$ with its image.
 The $\ZZ/5\ZZ$-, resp. $\alpha_5$-action on $Y'$ is induced by a linear
 action on $\PP^3$.
\end{Theorem}

\begin{Remark}
 If $X$ is a singular Godeaux surface, then $Y$ is smooth and 
 $Y'$, being the canonical model of $Y$, has at worst Du~Val singularities.
 We already noted in Remark \ref{singularremark} that $Y$ is not smooth
 if $X$ is supersingular.
\end{Remark}

\prf
Since $\Xcan$ has only rational singularities we obtain a Frobenius-equivariant
isomorphism $H^1(\OO_X)\iso H^1(\OO_{\Xcan})$.
Hence the associated torsor $Y\to X$ is the pullback
of a torsor of the same type over $\Xcan$.
We denote by $\nu:Y\to Y'$ the morphism induced by this pullback.
Being the pullback of an ample invertible sheaf
on $\Xcan$, namely $\omega_{\Xcan}$ possibly tensorised with an invertible 
sheaf that is torsion,
it follows that $\omega_{Y'}$ is ample.

In particular, $Y'$ is isomorphic to its canonical model.
We get $\omega_Y\iso\nu^\ast\omega_{Y'}$ and 
$\nu_\ast\omega_Y\iso\omega_{Y'}$
by flat base change and the corresponding properties of $\omega_X$.
Hence $Y'$ is the canonical model of $Y$ and the canonical map
of $Y$ factors over the canonical map of $Y'$.

We note that here, and in the sequel, we will freely use Ekedahl's results
\cite{ek} 
about canonical models of surfaces of general type in positive
characteristic.

We have $\chi(\OO_Y)=5$ and hence $h^0(\omega_Y)\geq4$, whereas 
Proposition \ref{noether} tells us $h^0(\omega_Y)\leq4$, whence
$h^0(\omega_Y)=4$ and $h^1(\OO_Y)=0$.
Let $\varphi_1$ be the canonical map of $Y$, which is a possibly rational
map to $\PP^3$, and whose image is a surface by Proposition \ref{noether}.

Let $X$ be a supersingular Godeaux surface.
We denote by $\rho:X^{(p)}\to Y$ the map
induced by $\pi:Y\to X$ and the Frobenius morphism.
As in the proof of Proposition \ref{noether}, we write the 
linear subsystem $\rho^\ast H^0(\omega_Y)$ of $|pK_{X^{(p)}}|$
as $|M|+V$.
Using that $K_X$ is nef, we get $M^2\leq (pK_X)^2=25$
and using Ekedahl's inequality \cite[Proposition 0.1.3]{ek}
we get
$M^2/(p\cdot \deg\varphi_1)\geq h^0(\omega_Y)-2=2$,
whence $\deg\varphi_1\leq2$.
Suppose we had $\deg\varphi_1=2$.
Then
\begin{equation}
 \label{estimates}
  25\,=\,(p K_X)^2\,=\,M^2\,+\,M\cdot V\,+\, pK_X\cdot V \,.
\end{equation}
First, we show $K_XV=0$ (still assuming $\deg\varphi_1=2$).
Suppose that $K_X V>0$.
Since pluricanonical divisors are connected, we have 
$MV>0$.
Furthermore we have $M^2\geq p\deg\varphi_1\deg(Z)\geq20$.
Plugging these estimates into (\ref{estimates}) we obtain
a contradiction.
Hence $K_XV=0$, which implies that $V$ is a sum 
of fundamental cycles.
Let $C$ be the unique effective divisor in $|K_X|$.
Since $F^\ast C$ lies in $\rho^\ast H^0(\omega_Y)$, the
fixed part $V$ is contained in $F^\ast C$.
Let $Z$ be a fundamental cycle in $V$.
Then the reduction of $\pi^\ast Z$ is contained in 
the base locus of $\omega_Y$.
But $\omega_Y\iso\pi^\ast\omega_X$, and hence
the {\it scheme} of base points of $H^0(\omega_Y)$
is $\alpha_p$-invariant.
Hence, $\pi^\ast Z$ is contained in this base scheme
from which we conclude that 
$\rho^\ast\pi^\ast Z=F^\ast Z=pZ$ is contained
in $V$.
Hence the fixed part is of the form $V = \sum_i n_iZ_i$
for disjoint fundamental cycles $Z_i$ and integers $n_i\geq p$.
We compute $0<M^2=(pK_X-V)^2=p^2 - 2\sum_i n_i^2$, which
is negative unless $V$ is empty.
Hence $V$ is empty, $M^2=p^2=25$, $\deg\varphi_1=2$ and
$\deg Z=2$.
In particular, $M$ has base points.

Let $\ideala:={\rm Bs}(\omega_Y)$ be the ideal sheaf of the base locus
of $H^0(\omega_Y)$.
As $M$ has base points there exists a closed point $y\in Y$ such that
$\ideala\subseteq\idealm_y$.
We already noted that $\ideala$ is $\alpha_p$-invariant, from which
we conclude that 
$\ideala\subseteq\pi^{-1}(\idealm_{\pi(y)})\cdot\OO_Y$.
Hence the ideal sheaf of the base locus of $\rho^\ast H^0(\omega_Y)$
is contained in $F^{-1}(\idealm_{\pi(y)})\cdot\OO_{X^{(p)}}$.
This means that $M$ has a base point of multiplicity $p^2$, which
is absurd as $M$ defines a morphism onto a surface and $M^2=p^2$.
This contradiction shows that $\varphi_1$ cannot be of degree $2$.

Hence $\varphi_1$ is birational.
Let $C$ be the unique effective divisor of $|K_X|$.
Since $K_X^2=1$, it follows that $C$ is of the form $C'+Z$,
where $C'$ is a reduced and irreducible curve and $Z$ is a sum
of fundamental cycles.
We write $\rho^\ast H^0(\omega_Y)=|M|+V$.
As before, the support of $V$ is contained in the support of $C$.
If $V$ contains $C'$ then $V$ also contains $pC'$ (using
$\alpha_p$-invariance as above) and $M$ is linearly equivalent
to a sum of fundamental cycles, which is absurd.
Hence $M$ is of the form $pK_X-\sum n_iZ_i$ for disjoint fundamental
cycles $Z_i$ and integers $n_i\geq p$, using again 
$\alpha_p$-invariance.
If we assume that the fixed part is not empty, we obtain the 
contradiction $M^2<0$.
Hence $V$ is empty and $M^2=p^2=25$.
As base points have to be counted with multiplicity $p^2$ and
the image of $\varphi_1$ is a surface, we conclude that $M$
has no base points.
Hence the image $W$ of $\varphi_1$ is a quintic hypersurface
in $\PP^3$.

By construction, we can identify $H^0(\omega_Y)$ 
and $H^0(\OO_W(1))$.
For every integer $n\geq1$, the map $\mu_n$ in the diagram
$$\begin{array}{ccc}
   H^0(W,\,\OO_W(1))^{\otimes n} &\iso& H^0(Y,\,\omega_Y)^{\otimes n}\\
   \downarrow {\mu_n} && \downarrow {\mu_n'}\\
  H^0(Z,\,\OO_W(n)) &\stackrel{\exists}{\dashrightarrow}&H^0(Y,\,\omega_Y^{\otimes n}) 
  \end{array}$$
is surjective, from which we conclude the existence of the dotted homomorphism
making the diagram commutative.
Since $\varphi_1$ is an isomorphism over an open dense set, 
we obtain an inclusion of $H^0(W,\OO_W(n))$ into $H^0(Y,\omega_Y^{\otimes n})$
for all $n\geq0$.
We have already seen $h^0(\omega_Y)=4$.
Now, $\omega_Y$ is the pullback of $\omega_X$ possibly tensorised with
an invertible sheaf that is torsion.
Hence, by the description of $\alpha_p$-torsors in
\cite[Proposition I.1.7]{ek}, we conclude that $\pi_\ast(\omega_Y^{\otimes n})$
is a successive extension of invertible sheaves on $X$, all of which
are numerically equivalent to $\omega_X^{\otimes n}$.
By Ekedahl's version of Mumford--Ramanujam vanishing
\cite[Proposition II.1.7]{ek}, the $H^1$ of these invertible sheaves
on $X$ is zero for all $n\geq2$.
From this we conclude that $h^1(\omega_Y^{\otimes n})=0$ for all
$n\geq2$ and applying Riemann--Roch we compute
$h^0(\omega_Y^{\otimes n})=5+5n(n-1)/2$ for $n\geq2$.
But these are the same dimensions as
$h^0(W,\OO_W(n))=h^0(W,\omega_W^{\otimes n})$.
Hence the injective homomorphism
$H^0(W,\OO_W(n)) \to H^0(Y,\omega_Y^{\otimes n})$
induced by $\varphi_1$ is an isomorphism for all $n\geq0$.
Hence $\varphi_1$ induces an isomorphism of the canonical models of $W$
and $Y$, i.e., an isomorphism of $Y'$ with $W$, as $W$ is isomorphic
to its canonical model.

Since $\mu_n$ and $\mu_n'$ are surjective, it follows that the $\alpha_p$-action
on $Y'$ and $W$ are completely determined by the linear action on 
$H^0(Y',\omega_{Y'})\iso H^0(W,\omega_W)$.
In particular, the $\alpha_p$-action on $W$ is induced by a linear action
on $\PP(H^0(W,\omega_W)^\vee)\iso\PP^3$.

If $X$ is a singular Godeaux surface, then $Y$ is a smooth surface.
Let $\varpi:\tilde{Y}\to Y$ be a blow-up such that the
movable part $|L|$ of $|\varpi^\ast K_Y|$ has no base points.
Since the image of $\varphi_1$ is a surface, a generic member 
$C$ of $L$ is a reduced and irreducible curve, hence
$h^0(C,\OO_C)=1$.
We also have $h^1(\OO_{\tilde{Y}})=0$, which then implies
$h^1(\OO_{\tilde{Y}}(-C))=0$.
Having thus established Kodaira vanishing for $L$
on $\tilde{Y}$ by hand, we can argue as in
the proof of \cite[Lemma 2]{horquint}
to conclude that $|K_Y|$ has no fixed part and at 
most one base point.

Since $\omega_Y$ is the pull-back of an invertible sheaf on $X$,
the group $\ZZ/5\ZZ$ acts on $H^0(Y,\omega_Y)$.
Hence, if $|K_Y|$ had a base point then the free 
$(\ZZ/5\ZZ)$-action would give at least $5$ base points.
Hence $|K_Y|$ has no base points and defines a morphism.
Its image $W$ is a surface, and since $K_Y^2=5$,
this map is birational onto a quintic hypersurface in $\PP^3$.
We leave the rest to the reader.
\qed\medskip

As already mentioned, existence of singular Godeaux surfaces has
been settled by Miranda \cite{mir}:

\begin{Theorem}[Miranda] 
  There do exist singular Godeaux surfaces in characteristic $5$.
 \qed
\end{Theorem}

It remains to establish the existence of supersingular Godeaux surfaces:

\begin{Theorem}
  \label{supersingexist}
  There do exist supersingular Godeaux surfaces in characteristic $5$.
\end{Theorem}

\prf
On $\PP^3$ we choose homogeneous coordinates
$x_0$, $x_1$, $x_2$, $x_3$ and consider the derivation
$$
\delta \,:=\, x_1\,\frac{\partial}{\partial x_0} 
\,+\,x_2\,\frac{\partial}{\partial x_1}
\,+\,x_3\,\frac{\partial}{\partial x_2}\,,
$$
which defines an additive vector field
in characteristic $5$, i.e., $\delta^{[5]}=0$.
This vector field has precisely one singular point, namely $[1:0:0:0]$.
Hence $\delta$ gives rise to an $\alpha_5$-action on $\PP^3$ with one fixed point.

Let $R:=k[x_0,...,x_3]$ be the homogeneous coordinate ring of $\PP^3$ and
let $R^\delta$ be the fixed ring of $\delta$.
Let $V$ be the vector space of elements of degree $5$ in $R^\delta$.
Since $x_i^5\in V$ for $i=0,...,3$, the rational map determined by $V$
$$
\begin{array}{ccccc}
 \varphi &:& \PP^3 &\dashrightarrow&\PP(V)
\end{array}
$$
is in fact a morphism.
We want to show that $\varphi$ can be identified with the quotient 
map $\PP^3\to\PP^3/\alpha_5$, at least outside the point
$[1:0:0:0]$.
To do so, we compute in local charts.

We set $x:=x_0/x_3$, $y:=x_1/x_3$ and $z:=x_2/x_3$
and define $A:=k[x,y,z]$.
Then $\delta$ becomes $y\partial_x\,+\,z\partial_y\,+\,\partial_z$ and we
consider 
$$
\begin{array}{lcl}
 \mbox{homogeneous of degree $5$} & \mbox{ \quad } & \mbox{inhomogeneous}\\
 F\,:=\,x_1x_3^4\,+\,2\,x_2^2x_3^3& &
 f\,:=\,y\,+\,2\,z^2\\
 G\,:=\,-x_0x_3^4\,+\,x_1x_2x_3^3\,+\,3\,x_2^3x_3^2& &
 g\,:=\,-x\,+\,y\,z\,+\,3\,z^3,
\end{array}
$$
which lie in the kernel of $\delta$.
An easy calculation shows that $B:=A^5[f,g]$ is smooth over $k$.
The fields of fractions of $A^\delta$ and $B$  are purely inseparable 
extensions of $k(x^5,y^5,z^5)$
of degree $5^2$ and height one.
On the other hand, $B$ is contained in $A^\delta$ and hence they
have the same field of fractions.
Moreover, since both rings are normal, they must be equal.
Hence, outside $\{x_3=0\}$ the map $\varphi$ induces an isomorphism
of its image with $\PP^3/\alpha_5$.

Now, we set $\tilde{x}:=x_0/x_2$, $\tilde{y}:=x_1/x_2$ and 
$\tilde{z}:=x_3/x_2$ and define
$\tilde{A}:=k[\tilde{x},\tilde{y},\tilde{z}]$.
Then $\delta$ becomes 
$(\tilde{y}-\tilde{x}\tilde{z})\partial_{\tilde{x}}+
(1-\tilde{y}\tilde{z})\partial_{\tilde{y}}-
\tilde{z}^2\partial_{\tilde{z}}$.
Let $\tilde{B}$ be the $\tilde{A}^5$-algebra generated by the dehomogenised
elements of $V$.
We consider $\Spec\tilde{A}^\delta\to\Spec\tilde{B}\to\Spec\tilde{A}^5$.
From the computations above we infer that 
for every closed point of $\Spec\tilde{A}^5$ whose 
$\tilde{z}^5$-coordinate is not equal to zero, there exists an 
open neighbourhood over which $\Spec\tilde{B}$ and $\Spec\tilde{A}^\delta$ are isomorphic.
Next, we consider the elements
$$
\begin{array}{lcl}
 \mbox{homogeneous of degree $5$} & \mbox{ \quad } & \mbox{inhomogeneous}\\
 x_3^{-5}F^2 
 & &
 a\,:=\,-\tilde{z}\,-\,\tilde{y}\tilde{z}^2\,+\,\tilde{y}^2\tilde{z}^3 \\

 x_3^{-10}(F^2G-x_2^5x_3^5F)
 & &
 b\,:=\,\tilde{x}\,+\,2\tilde{y}^2\,+\,\tilde{y}^3\tilde{z}\,+\,
 \tilde{x}\tilde{y}\tilde{z}\,-\,\tilde{x}\tilde{y}^2\tilde{z}^2
\end{array}
$$
which lie in the kernel of $\delta$.
We define $\tilde{C}:=\tilde{A}^5[a,b]$.
Then we have an inclusion $\tilde{C}\subseteq\tilde{B}\subseteq\tilde{A}^\delta$
and all three rings have the same field of fractions.
A straight forward computation shows that
for every closed point of $\Spec\tilde{A}^5$ with $\tilde{z}^5$-coordinate equal to zero,
there exists an open neighbourhood, s.th. all closed points of $\Spec\tilde{C}$ lying
above this neighbourhood are smooth over $\Spec k$.
Hence over this neighbourhood $\Spec\tilde{A}^\delta$ and $\Spec\tilde{C}$ are
isomorphic, and in particular $\Spec\tilde{A}^\delta$ and $\Spec\tilde{B}$ are isomorphic
over this neighbourhood.
We conclude that $\tilde{B}$ and $\tilde{A}^\delta$ are isomorphic
and that $\varphi$ can be identified with the quotient map 
by $\alpha_5$ outside $\{x_2=0\}$.

Next, we set $\hat{x}:=x_0/x_1$, $\hat{y}:=x_2/x_1$ and 
$\hat{z}:=x_3/x_1$ and define
$\hat{A}:=k[\hat{x},\hat{y},\hat{z}]$.
Then $\delta$ becomes 
$(1-\hat{x}\hat{y})\partial_{\hat{x}}+
(\hat{z}-\hat{y}^2)\partial_{\hat{y}}-
\hat{y}\hat{z}\partial_{\hat{z}}$.
We let $\hat{B}$ be the $\hat{A}^5$-algebra generated by the dehomogenised
elements of $V$.
We want to show that $\Spec\hat{A}^\delta$ and $\Spec\hat{B}$ are isomorphic,
which - by the previous computations - is clear over the open set of points 
lying above points of
$\Spec\hat{A}^5$ with $\hat{y}^5$-coordinate or $\hat{z}^5$-coordinate
not equal to zero.
By the previous considerations it is enough to find a
$\hat{A}^5$-subalgebra $\hat{C}$ of $\hat{B}$ with the 
same field of fractions such that $\Spec\hat{C}$ is smooth over 
$\Spec k$ at all closed points lying above points 
of $\Spec\hat{A}^5$ with
$\hat{y}^5$- and $\hat{z}^5$-coordinate equal to zero.
We consider
$$
\begin{array}{lcl}
 \mbox{homogeneous of degree $5$} & \mbox{ \quad } & \mbox{inhomogeneous}\\
  x_3^{-15}(F^4+2x_3^5FG^2+x_2^5x_3^{10}G)
 & & 
 c\,:=\,\hat{z}\,+\,2\hat{x}^2\hat{z}^2\,+\,\hat{x}\hat{y}\hat{z}\,-\,\hat{x}^2\hat{y}^2\hat{z}\\

 x_3^{-15}(G^3x_3^5+F^3G+F^2x_2^5x_3^5)
 & &
 d\,:=\,\hat{y}\,-\,\hat{x}\hat{z}\,-\,\hat{x}^2\hat{y}^3\,+\,3\hat{x}^2\hat{y}\hat{z}\,-\,
 \hat{x}^3\hat{z}^2\,+\,\hat{x}\hat{y}^2\\
\end{array}
$$
which lie in the kernel of $\delta$.
An easy computation shows that $\hat{C}:=\hat{A}^5[c,d]$ 
has the desired properties.
As before, we conclude that $\varphi$ can be identified with the quotient map 
by $\alpha_5$ outside $\{x_1=0\}$.

Hence $Z:={\rm im}(\varphi)$ coincides with $\PP^3/\alpha_5$ except possibly
at $\varphi(P)$, where $P:=[1:0:0:0]$.
We have also seen in the computations above that $Z$ is even smooth outside $\varphi(P)$.

Let $X$ be a generic hyperplane section of $Z$, which is a smooth
surface since $Z$ has only one isolated singular point.
The inverse image $Q$ of $X$ under $\varphi$ in $\PP^3$ 
is a $\delta$-invariant quintic and gives rise to
an $\alpha_5$-torsor $\pi:Q\to X$.

Being a quintic surface in $\PP^3$, we compute $\chi(\OO_Q)=K_Q^2=5$, where
the self-intersection of the canonical divisor is understood in the
sense of formula (\ref{selfintersection}).
Using (\ref{standard}) we obtain $\chi(\OO_X)=K_X^2=1$.
Since $\omega_Q=\pi^\ast(\omega_X)$ is ample, $X$ is a minimal surface 
of general type, hence a numerical Godeaux surface.
The existence of the $\alpha_5$-torsor $\pi$ shows that $X$ is supersingular.
\qed

\end{document}